\documentclass[11pt]{article}
\usepackage{xparse}
\usepackage{amsmath}
\usepackage{amssymb}
\usepackage{textcomp}
\usepackage[document]{ragged2e}
\usepackage{setspace}
\usepackage{comment}
\usepackage{slashed}
\usepackage{etoolbox}
\usepackage[affil-it]{authblk}

\usepackage[dvipsnames]{xcolor}

\usepackage[pdfencoding=auto, psdextra]{hyperref}
 \hypersetup{
     colorlinks=true,
     linkcolor=blue,
     citecolor = blue,      
     urlcolor=blue,
     }
\newenvironment{changemargin}[2]{
\begin{list}{}{
\setlength{\topsep}{0pt}
\setlength{\leftmargin}{#1}
\setlength{\rightmargin}{#2}
\setlength{\listparindent}{\parindent}
\setlength{\itemindent}{\parindent}
\setlength{\parsep}{\parskip}
}
\item[]}{\end{list}}

\newcommand{\bb}[1]{\mathbb{#1}}
\newcommand{\alg}[1]{\mathfrak{#1}}
\newcommand{\spin}{\mathfrak{spin}}
\newcommand{\su}{\mathfrak{su}}
\newcommand{\so}{\mathfrak{so}}
\newcommand{\grp}[1]{\operatorname{#1}}

\title{Sporadic Isogenies to the Quaternionic\newline Orthogonal Groups $\operatorname{SO}^*(2n)$}
\author{Craig M$^{\mathrm{c}}$Rae\thanks{Electronic address: \texttt{mcraec3@myumanitoba.ca}}}
\affil{University of Manitoba \\ Winnipeg, MB}

\begin{document}
\maketitle
  \pagenumbering{gobble}
  \pagenumbering{arabic} 
\setlength{\parindent}{1.5em}
\setlength{\parskip}{0.5ex}
\begin{changemargin}{-1cm}{-1cm}
\begin{abstract}
\normalsize\noindent For small dimensional Lie algebra’s there are many so-called accidental isomorphisms which give rise to double covers of special orthogonal groups --- Spin groups --- which happen to coincide with groups already belonging to another classification. The well known catalog of these sporadic isogenies given by Dr. Paul Garrett \cite{isogenies} keeps track of these facts for the complex orthogonal groups $\operatorname{SO}(n, \mathbb{C})$, and those real forms which are \textit{manifestly} real: $\operatorname{SO}(p,q)$, while avoiding the quaternionic real forms, $\operatorname{SO}^*(2n)$. This article serves to complement and complete the existing catalog by presenting the sporadic isogenies to the first four quaternionic orthogonal groups in one place, with contemporary proofs of those group and algebra homomorphisms. A review of the definition of $\operatorname{SO}^*(2n)$ is presented in a modern notation, emphasizing the concept of `quaternion reversion' as a useful, albeit redundant, second conjugation upon quaternions, helping to explicate the self-conjugacy of quaternionic representations. Lastly, the triality of $\operatorname{SO}^*(8)$ is explored in a similar manner to work done by the author in \cite{Triality}.
\end{abstract}
\end{changemargin}
\newpage 
\begin{changemargin}{-2cm}{-2cm}
\tableofcontents

\newpage

\baselineskip=1.2\baselineskip
\section*{Preamble}
\label{sec:Preamble}
\addcontentsline{toc}{section}{\nameref{sec:Preamble}}
For those looking for a quick reference, below is a table identifying algebra homomorphisms to $\alg{so}^*(2n)$, the unique quaternionic real form of $\alg{so}(2n,\mathbb{C})$, alongside the sporadic isogenies to the defining representations of $\grp{SO}^*(2n)$ in low dimension. These results are not novel, but the groups themselves not particularly well known, let alone the proofs of these homomorphisms; this article hopes to remedy this. 

\vspace{0.5em}
\begin{table}[h]
    \centering \large
        \begin{tabular}{|c|c|c|}
        \hline
         $n$ & $\alg{so}^*(2n)$ & $\grp{SO}^*(2n)$\\ \hline
         $1$ & $\so(2,\mathbb{R})$ & $\grp{SO}(2, \mathbb{R})$\\
         $2$ & $\su(2) \oplus \alg{sl}(2,\mathbb{R})$ & $\left(\grp{SU}(2) \times \grp{SL}(2,\mathbb{R}) \right)/\mathbb{Z}_2$ \\
         $3$ & $\alg{su}(3,1)$ & $\grp{SU}(3,1) / \mathbb{Z}_2$\\
         $4$ & $\alg{spin}(2,6)$ & $\grp{Spin}(2,6) /\mathbb{Z}_2 \>\cong \> \grp{SemiSpin}(2,6) $\\\hline
    \end{tabular}
    \begin{changemargin}{-1.1cm}{-1.1cm}
    \caption{The first four $\alg{so}^*(2n)$ algebras given in terms of the algebra of their double covering; i.e. $\spin^*(2n)$. To the authors knowledge, it is only these four cases where $\so^*(2n)$ coincides with another known Lie algebra. It is interesting to note that each of these Lie groups has an infinite fundamental group.}
    \label{tab:orth}
    \end{changemargin}
\end{table}
\begin{table}[h]
\centering \large
    \begin{tabular}{|c|ccc|}
        \hline Name & dim $d$ & rank & \# Compact Generators \\\hline
        $\alg{so}^*(2n)$ & $n(2n-1)$ & $n$ & $n^2$ \\\hline
    \end{tabular}
    \begin{changemargin}{-1.1cm}{-1.1cm}
    \caption{Basic facts about the Lie algebra's $\so^*(2n)$, the unique quaternionic real form of $\so(2n, \mathbb{C})$.}
    \label{tab:sostar}
    \end{changemargin}
\end{table}

\vspace{-0.5em}
\noindent
The quaternionic orthogonal groups in the defining representation are those quaternion valued matrices which preserve a \textit{reversion} symmetric product on $\mathbb{H}^n$:
\begin{equation}
\begin{split}
\grp{SO}^*(2n) &= \{O \in \grp{SL}(n, \mathbb{H}) \>\> | \>\> \widetilde{O}^\intercal \mathbb{I}_n O = \mathbb{I}_{n}\},\\
\alg{so}^*(2n) &= \{\>a \in \alg{sl}(n, \mathbb{H}) \>\> \>\> | \quad -\widetilde{a}^\intercal =  a \}.
\end{split}
\end{equation}
Reversion, denoted with a tilde $\widetilde{\>\>\>}$, is conjugation of only the $j$ component of a quaternion. Quaternion reversion gives a convenient way of talking about quaternion sesquilinear-products which are skew-Hermitian. Reversion is an anti-homomorphism on quaternions, and reversion-transpose is an anti-homomorphism on quaternion valued matrices. For a detailed exposé and motivation for the introduction of reversion, such as the relationship between reversion on $\mathbb{H}^n$ and quaternionic and orthogonal structures on $\mathbb{C}^{2n}$, see Sec.~(\ref{sec:Reversion}), Sec.~(\ref{sec:quatorthgrps}), and Sec.~(\ref{sec:waxdiff}). 

\noindent
The following four sections explicitly construct these algebras, and demonstrate these group and algebra homomorphisms --- they rely on a strong grasp of quaternionic representation theory. A review of these topics is given in Appendix~\ref{app:A}. The work throughout will make regular use of the homomorphism from $\grp{M}_n(\mathbb{H}) \cong \grp{M}_{2n}(\mathbb{C})$, given via element wise identification:
\begin{equation}
\label{eq:quat2Pauli}
    t+xi+yj+zk \quad \mapsto \quad \begin{pmatrix}
        t + zi & ix -y\\
        ix +y & t - zi
    \end{pmatrix}.
\end{equation}

\noindent\rule{16.5cm}{0.4pt}
\section{\texorpdfstring{$\grp{SO}^*(2) \cong \grp{SO(2, \mathbb{R})}$}{SO*(2) = SO(2, R)}}
In this case we are looking for $1\times 1$ quaternions $R$, satisfying $\widetilde{R}^\intercal R = 1$. The Lie algebra is then the set of quaternions $a$ satisfying $\widetilde{a} = -a$: this give a one dimensional Lie algebra, its elements all of the form $\theta j$. Clearly the exponential map sends us to the unit circle in the $1-j$ plane of quaternions, and so the group is isomorphic to $\grp{U}(1)$. Equivalently, mapping $j$ to its standard complex matrix via Eq.~(\ref{eq:quat2Pauli}) we have that the group elements are 
\begin{equation}
    \grp{SO}^*(2) \ni e^{\theta j} \leftrightarrow\exp\begin{pmatrix}
        0 & -\theta \\
        \theta & 0
    \end{pmatrix} = \begin{pmatrix}
        \cos\theta & -\sin \theta\\
        \sin\theta & \cos\theta
    \end{pmatrix} \in \grp{SO}(2, \mathbb{R}).
\end{equation}
Thus $\grp{SO}^*(2) \cong \grp{SO(2, \mathbb{R})} \cong \grp{U}(1)$.

\noindent\rule{16.5cm}{0.4pt}
\section{\texorpdfstring{$\grp{SO}^*(4) \cong \left(\grp{SU}(2) \times \grp{SL}(2, \mathbb{R}) \right)/\mathbb{Z}_2$}{SO*(4) = (SU(2) x SL(2, R)) / Z2}}
This group is $6$ dimensional, rank $2$ and its largest compact subgroup is four dimensional. Focusing firstly upon the Lie algebra, we have that
\begin{equation}
\alg{so}^*(4) = \{A \in \alg{sl}(\mathbb{H}^2) \> \> | \>\> -\widetilde{A}^\intercal = A \}.
\end{equation}
One choice of basis for this may be given by elements $A_i$:
\begin{equation}
\sum_{i=1}^6\theta_i A_i = \frac{1}{2} \begin{pmatrix}
     (\theta_2+\theta_5) j & \theta_1 +\theta_4 i-\theta_3 j+\theta_6 k \\
    -\theta_1-\theta_4 i -\theta_3 j-\theta_6 k & (-\theta_2+\theta_5) j \\
\end{pmatrix}.
\end{equation}
In this basis the generators $A_1,A_2,$ and $A_3$, separate from $A_4, A_5$, and $A_6$: these two sets of three are commuting sub-algebra's. In particular it is easy to verify that $A_1,A_2,$ and $A_3$ form an $\su(2)$ sub-algebra, while $A_4,A_5,$ and $A_6$ form an $\alg{sl}(2,\mathbb{R})$ sub-algebra. Therefore at the level of the Lie algebra, it is plain to see 
\begin{equation}
    \alg{so}^*(4) = \su(2)\oplus \alg{sl}(2, \mathbb{R}) .
\end{equation}
However the situation at the level of the group is not so clean. Even though the algebra is reducible to a direct sum of simple algebra's, by explicit computation one can show the given representation of $\grp{SO}^*(4)$, is irreducible: when promoted to the corresponding representation over $\mathbb{C}$ instead of $\mathbb{H}$ (again via Eq.~(\ref{eq:quat2Pauli})), there are only trivial intertwiners and so there are no invariant subspaces of the group action. 

This is exactly analogous to what happens with the vector representation of the compact real form of this group, $\grp{SO}(4, \mathbb{R})$, which is famously isomorphic to $\left(\grp{SU}(2)\times \grp{SU}(2)\right) /\mathbb{Z}_2$. The story is the same here. If we consider the reducible representation given by the action of $\grp{SU}(2) \oplus \grp{SL}(2,\mathbb{R})$ acting on $\mathbb{C}^2 \oplus \mathbb{C}^2 \cong \mathbb{C}^4$, this representation has a center of $\mathbb{Z}_2 \times \mathbb{Z}_2$, one factor from each component. This reducible representation may be written in the following basis $S_i$:
\begin{equation}
    \sum_{i=1}^{6} \theta_i S_i= \frac{1}{2}\begin{pmatrix}
     i \theta _3 & -\theta _2+i \theta _1 & 0 & 0 \\
     \theta _2+i \theta _1 & -i \theta _3 & 0 & 0 \\
     0 & 0 & -\theta _6 & -\theta _4-\theta _5 \\
     0 & 0 & \theta _5-\theta _4 & \theta _6 \\
    \end{pmatrix}.
\end{equation}
The top left are the standard generators of $\su(2)$, and the bottom right, $\alg{sl}(2,\bb{R})$. Importantly the bases $A_i$ and $S_i$ are constructed to have identical Lie bracket, and so mapping between them is an algebra automorphism. Define this algebra automorphism $\sigma: S \mapsto A$. Let us now see what the `difference' is between these group representations via the exponential map from the algebra: i.e. what is the kernel of the induced homomorphism by $\sigma$ upon the group representations? Inspecting some rotations by $2 \pi$ we have:
\begin{equation}
\begin{split}
    U_1 = \exp\left( 2 \pi S_1 \right) &= \begin{pmatrix}
        -\mathbb{I}_2 & 0 \\ 
        0 & \mathbb{I}_2
    \end{pmatrix}, \quad \stackrel{\Large\sigma}{\mapsto}\quad  R_1 =\exp\left( 2 \pi A_1\right) = \begin{pmatrix}
        -\mathbb{I}_2 & 0 \\ 
        0 & -\mathbb{I}_2
    \end{pmatrix}, \\
    U_5=   \exp\left( 2 \pi S_5 \right) &= \begin{pmatrix}
        \mathbb{I}_2 & 0 \\ 
        0 & -\mathbb{I}_2
    \end{pmatrix}, \quad \stackrel{\Large\sigma}{\mapsto} \quad  R_5 = \exp\left( 2 \pi A_5 \right) = \begin{pmatrix}
        -\mathbb{I}_2 & 0 \\ 
        0 & -\mathbb{I}_2
    \end{pmatrix},
\end{split}
\end{equation}
In particular the induced group homomorphism is not an automorphism: 
\begin{equation}
    U_1 U_5 = -\mathbb{I}_4, \quad \stackrel{\Large\sigma}{\mapsto } \quad R_1 R_5 = \mathbb{I}_4.
\end{equation}
This shows the algebra automorphism, mapping us between these representations corresponds to a double covering at the level of the group. I.e. for every element of $\grp{SO}^*(4)$, there are two corresponding elements of $\grp{SU}(2) \times \grp{SL}(2,\mathbb{R})$, related by application of an element in the center. In particular there is a group isomorphism
\begin{equation}
    \grp{SO}^*(4) \cong \left(\grp{SU}(2) \times \grp{SL}(2,\mathbb{R})\right)/\mathbb{Z}_2,
\end{equation}
where the modulo by $\mathbb{Z}_2$ identifies the every matrix with its negation. With this understood we may identify the largest compact subgroup as the group $\grp{S}(\grp{U}(2)\times \grp{SO}(2)) \cong \grp{S}(\grp{U}(2)\times \grp{U}(1))$, which is interesting in as much as it relates to the physics of flavourdynamics, and gives a somewhat natural place for the mixing of a $\grp{U}(1)$ and $\grp{SU}(2)$ symmetry to arise from a simple Lie group. This pattern will come up again in a far more compelling way in the sext section. 

\noindent\rule{16.5cm}{0.4pt}
\section{\texorpdfstring{$\grp{SO}^*(6) \cong \grp{SU}(3,1) /\mathbb{Z}_2$}{SO*(6) = SU(3,1) / Z2}}
\subsection*{Basis for the Lie algebra \texorpdfstring{$\alg{su}(3,1)$}{su(3,1)}}
The group $\grp{SU}(3,1)$ in the defining representation is given by
\begin{equation}
    \grp{SU}(3,1) = \{S \in \grp{SL(4,\mathbb{C})} \> \> | \>\> S^\dagger \mathbb{I}_{3,1} S = \mathbb{I}_{3,1}\}.
\end{equation}
Here $\mathbb{I}_{p,q}$ is the diagonal matrix with $p$ entries being $1$ and $q$ entries $-1$. The group is $15$ dimensional, rank $3$ and its largest compact subgroup is $\grp{S}(\grp{U}(1)\times \grp{U}(3))$ which is $9$ dimensional. The Lie algebra $\alg{su}(3,1)$ is then given by those matrices $s$ satisfying:
\begin{equation}
    \alg{su}(3,1) = \{s \in \alg{sl}(4,\mathbb{C}) \> \> | \>\> \mathbb{I}_{3,1} s \>\mathbb{I}_{3,1} = -s^\dagger  \},
\end{equation}
for which an orthogonal basis of fifteen 4 by 4 generators, $s_i$, may be given by:
\begin{equation}
    s_i = \begin{pmatrix}
        \lambda_i & 0 \\
        0 & 0
    \end{pmatrix}, \qquad i \in (1,8), \qquad s_{15} = \frac{i }{2\sqrt{6}}\begin{pmatrix}
  \mathbb{I}_3 & 0 \\
 0 & - 3 \\
\end{pmatrix},
\end{equation}
\begin{equation}
\sum_{i=9}^{14} s_{i} \theta_i = \frac{1}{2}\begin{pmatrix}
 0 & 0 & 0 & -\theta _{13}+i \theta _{14} \\
 0 & 0 & 0 & \theta _{11}-i \theta _{12} \\
 0 & 0 & 0 & -\theta _9+i \theta _{10} \\
 -\theta _{13}-i \theta _{14} & \theta _{11}+i \theta _{12} & -\theta _9-i \theta _{10} & 0 \\
\end{pmatrix}.
\end{equation}
Here $\lambda_i$ are the mathematicians Gell-Mann matrices, given in appendix~(\ref{sec:su3}). It is fun to note that the manifestly real sub-algebra here is precisely the vector representation of $\so(3,1)$. The generators given are normalized to $\grp{Tr}(s_i s_i) = \pm\frac{1}{2} $.

\subsection*{Basis for the Lie algebra \texorpdfstring{$\alg{so}^*(6)$}{so*(6)}}
In the defining representation on $\mathbb{H}^3$, $\grp{SO}^*(6)$ is the group of quaternion-linear transformations satisfying:
\begin{equation}
\{A \in \grp{SL}(\mathbb{H}^3)\> \>| \>\>   \widetilde{A}^\intercal \mathbb{I}_{3} A = \mathbb{I}_{3} \}.
\end{equation}
The group is $15$ dimensional, rank $3$ and its largest compact subgroup is $9$ dimensional. The Lie algebra $\alg{so}^*(6)$ is then given by those matrices $a$ satisfying:
\begin{equation}
\alg{so}^*(6) = \{a \in \alg{sl}(\mathbb{H}^3) \> \> | \>\> -\widetilde{a}^\intercal = a  \},
\end{equation}
for which an (admittedly messy looking) orthogonal basis of fifteen 3 by 3 generators, $a_i$, may be given by:
\begin{equation*}
\begin{split}
\sum_{i=1}^{15} a_i \theta_i = \frac{1}{4}\begin{pmatrix}
        j\left(\sqrt{2} \>\theta_{15}-2\theta_8\right)/\sqrt{3}  & \theta_7 + \theta_{13} i + \theta_6 j + \theta_{14} k & \theta_5 + \theta_{11} i + \theta_4 j + \theta_{12} k \\
        -\theta_7 - \theta_{13} i + \theta_6 j - \theta_{14} k & j\left(\sqrt{2} \>\theta_{15}+\theta_8- \sqrt{3}\theta_3\right)/\sqrt{3} & \theta_{2} + \theta_{9} i + \theta_{1} j + \theta_{10} k\\
        -\theta_5 -\theta_{11} i + \theta_4 j - \theta_{12} k  & -\theta_{2} -\theta_{9} i + \theta_{1} j - \theta_{10} k & j\left(\sqrt{2} \>\theta_{15}+\theta_8+ \sqrt{3}\theta_3\right)/\sqrt{3}
    \end{pmatrix}.
\end{split}
\end{equation*}
The basis given is normalized such that $\grp{Tr}(a_i a_i) = \pm 1$. Let us make the construction more palatable by once again mapping to the complex case via Eq.~(\ref{eq:quat2Pauli}), embedding the given quaternionic basis into $\mathfrak{so}(6,\mathbb{C})$ (abusing notation and keeping the name of the basis $a$). In this setting, as usual one will find the group elements satisfy simultaneously an orthogonal, quaternionic, and Hermitian structure. In an effort to make the group `appear' more unitary we will diagonalize the Hermitian structure, employing a change of basis $U$ (given explicitly in Eq.~(\ref{changeofbasis}) in the appendix) on the Lie algebra representation. Our orthonormal Lie algebra generators $a_i$ now take the following convenient form:
\begin{equation}
\begin{split}
    a_i = \begin{pmatrix}
        \lambda_i & 0 \\
        0 & \lambda_i^*
    \end{pmatrix}&, \qquad i \in (1,8), \qquad a_{15} = \frac{i}{\sqrt{6}} \begin{pmatrix}
        \mathbb{I}_3 & 0 \\
        0 & -\mathbb{I}_3
    \end{pmatrix}, \\
    a_{9} = \begin{pmatrix}
        0 & -\lambda_2 \\
        \lambda_2 & 0 \\
    \end{pmatrix}, \qquad a_{11} &= \begin{pmatrix}
        0 & -\lambda_5 \\
        \lambda_5 & 0 \\
    \end{pmatrix}, \qquad a_{13} = \begin{pmatrix}
        0 & -\lambda_7 \\
        \lambda_7 & 0 \\
    \end{pmatrix}, \\
    a_{10} = \begin{pmatrix}
        0 & i\lambda_2 \\
        i\lambda_2 & 0 \\
    \end{pmatrix}, \qquad a_{12} &= \begin{pmatrix}
        0 & i\lambda_5 \\
        i\lambda_5 & 0 \\
    \end{pmatrix}, \qquad a_{14} = \begin{pmatrix}
        0 & i\lambda_7 \\
        i\lambda_7 & 0 \\
    \end{pmatrix}.
\end{split}
\end{equation}
Again the $\lambda_i$ are the mathematicians Gell-Mann matrices. It is clear the first $8$ elements generate a reducible representation of $\grp{SU(3)}$: $\mathbf{3}\oplus\overline{\mathbf{3}}$, and $a_{15}$ commutes with these giving us a $\grp{U}(3)$ in either block, though the phases in either block must be inverse to preserve the overall unit determinant. This explicitly gives the generators of the $9$ dimensional maximally compact subgroup which is precisely $\grp{S}(\grp{U}(1)\times \grp{U}(3))$. The remaining $6$ generators are our `boosts', as they are non-compact and mix up the nice subspaces of our compact subgroup, exactly analogous to how boosts mix up the compact $\grp{SO}(3)$ subspaces of the Faraday tensor (the adjoint representation) in $\grp{SO}(3,1)$. In fact the manifestly real sub-algebra here is precisely $\alg{so}(3,1)$ in the adjoint representation. 

\subsection*{Double Covering}
The relationship between the above two algebras is that they are identical. It is not a coincidence that they share rank, dimension, and maximal compact subgroups. One can see that the mapping $s_i \mapsto a_i$ is a Lie algebra automorphism: in the bases constructed the commutation relations are identical. That is to say, with:
\begin{equation}
     \left[s_i,s_j\right] = f_{ij}^{\>\>\>k}  s_k, \quad \left[a_i,a_j\right] = c_{ij}^{\>\>\>k} a_k, \quad \textrm{then}\quad  f_{ij}^{\>\>\>k}  = c_{ij}^{\>\>\>k} \quad \forall i,j,k.
\end{equation}The representation given of $\grp{SO}^*(6)$ can be shown to be precisely $\wedge^2 \grp{SU}(3,1)$. In particular one can see it is a $2:1$ covering from the group centers. The determinant constraint of an element of the center of $\grp{SU}(3,1)$ demands the center is at most the four roots of unity. We can see
\begin{equation}
    \grp{exp}\left(s_{15} \> \sqrt{6} \>\pi \right) = i \>\mathbb{I}_4,
\end{equation}
and so the center of the group is $\mathcal{Z}(\grp{SU}(3,1)) = \mathbb{Z}_4$, generated by the above element. However the same algebra element mapped into the quaternionic orthogonal group gives
\begin{equation}
    \grp{exp}\left(a_{15} \> \sqrt{6} \>\pi \right) = -\>\mathbb{I}_6,
\end{equation}
and so $\mathcal{Z}(\grp{SO}^*(6)) = \mathbb{Z}_2$, the covering map sending each element of the center to its square. Thus we have an isogeny showing $\grp{SU}(3,1)$ is the double cover of $\grp{SO}^*(6)$. One could equivalently say $\grp{Spin}^*(6) = \grp{SU}(3,1)$. 

The group $\grp{SU}(3,1)$ has a surprisingly similar structure to the basic algebras in fundamental physics. As noted above the real part of the algebra can be chosen to be isomorphic to the standard representations of the Lorentz group, naturally giving rise to vectors and field strength tensors. Meanwhile the maximal compact subgroup gives precisely the structures of the electromagnetic and chromodynamic forces, and in such a way that they must to some extent `know about' each other to preserve the unit determinant, giving a plausible cause to why colour singlets know anything about the charge of the electron. In work such as \cite{Wilson2020}, the merits and obstructions to an $\grp{SL}(4, \mathbb{R})$ unification scheme have been discussed, and as $\grp{SU}(3,1)$ is another real form of the same group, it may be worth while in future work to consider similar arguments for it. 

\noindent\rule{16.5cm}{0.4pt}
\section{\texorpdfstring{$\grp{SO}^*(8) \>\cong \>\grp{SemiSpin}(2,6) $}{SO*(8) = SemiSpin(2,6)}}
This section follows the conventions, notations, and definitions used in \cite{Triality}, a related article by the author. To inspect the stated isomorphism, first we can build the Clifford algebra $C\ell(2,6)$, its generators called $\Gamma_i$, to build a Lie algebra basis of $\spin(2,6)$, and then see it is `quaternion orthogonal'. We can build this conveniently out of a Clifford algebra basis $g_i$ for $C\ell(7,0)$:\footnote{The matrices can be found explicitly in Appendix \ref{app:cl7}.}

\begin{equation}
\label{eq:CL26}
\Gamma_0 = \begin{pmatrix}
0 & \mathbb{I}_8 \\
\mathbb{I}_8 & 0
\end{pmatrix}, \quad \Gamma_{1} = i\begin{pmatrix}
0 & g_1 \\
-g_1 & 0
\end{pmatrix},\quad \Gamma_{\mu} =  \begin{pmatrix}
0 & g_\mu \\
-g_\mu & 0
\end{pmatrix} \quad 2 \leq \mu \leq 7.
\end{equation}
The matrices mutually anti-commute and square to $+1,+1,$ and then $-1$ respectively. By the standard theory of Clifford algebras \cite{CliffAlg}, a basis for the $\spin(2,6)$ algebra is then given by
\begin{equation}
\label{LRdef}
S_{ij} = \frac{1}{2}\Gamma_i \Gamma_j, \quad 0 \leq i < j \leq 7, \qquad \grp{span}\{S_{ij}\} = \spin(2,6).
\end{equation}
We will also flip the signs of the following generators for convenience: 
\begin{equation}
\begin{split}
    S_{12} &\mapsto -S_{12}, \quad S_{15} \mapsto -S_{15}, \quad S_{17} \mapsto -S_{17}, \\
    S_{24} &\mapsto -S_{24}, \quad S_{26} \mapsto -S_{26}, \quad S_{47} \mapsto -S_{47}. \\
\end{split}
\end{equation}

\noindent
The algebra is $28$ dimensional, rank $4$ and its largest compact sub-algebra is that of $\grp{S}(\grp{SO}(2)\times \grp{SO}(6)) \cong \grp{S}(\grp{U}(1)\times \grp{U}(4))/\mathbb{Z}_2$, which is $16$ dimensional.\footnote{While less obviously connected to physics, the structure could still give an interesting unification scheme, as all the required sub-groups are present, alongside triality.} Note that as the matrices in Eq.~(\ref{eq:CL26}) are block anti-diagonal, their products will be block diagonal and so this $16$ dimensional representation of the Lie algebra is reducible and the matrices $S_{ij}$ may be decomposed into these blocks. We will call the top left block the left handed block, its generators denoted $L_{ij}$, and the bottom right block the right handed block, and its generators denoted $R_{ij}$. These left and right handed representations of $\spin(2,6)$ are distinct irreducible representations.\footnote{The corresponding group representations will be the so called SemiSpin representations, due to the fact that the center of the spin group is $\mathbb{Z}_2 \times \mathbb{Z}_2$, and so no irrep will be faithful, always missing at least half the center. The reducible representation however, given by the $S_{ij}$ would give a faithful representation of the spin group.} Physically one may think of these representations as being reflections of one another, as the inequivalent left and right handed representations swap when one sends $\Gamma_0 \mapsto -\Gamma_0$, which is a reflection on the $8$ dimensional vector space the Clifford algebra is initially defined on. It can be seen that both of these representations are self-dual, and self-conjugate. 

Focusing in particular on the left handed representation, let us write an arbitrary element of the Lie algebra in this representation as $s(\vec{\theta}) = \sum \theta_{ij} L_{ij}$, where $\theta_{ij}$ are the boost and rotation parameters for the $ij$-plane. Then with
\begin{equation}
    \varepsilon = \begin{pmatrix} 0 & -1 \\
    1 & 0\end{pmatrix}, \qquad J = \operatorname{diag}\{\varepsilon,\varepsilon,\varepsilon,\varepsilon\} ,\qquad g = \mathbb{I}_8,
\end{equation}
one can show that this arbitrary element $s$ in this representation of the Lie algebra satisfies the following properties:
\begin{equation}
    J s^* J^{-1}= s , \qquad g \left(-s^\intercal \right) g^{-1}= s.
\end{equation}
The first property tells us the Lie algebra representation is self-conjugate and so must be equivalent to a Lie algebra over the reals or quaternions; and  since $J J^* = -1$ we can conclude the representation is quaternionic. The second property tells us the generators are anti-symmetric and so are a sub-algebra of $\alg{so}(8, \mathbb{C})$. Thus we have simultaneous orthogonal and quaternionic structures, and can conclude we must be working with a representation of $\so^*(8)$.\footnote{Technically as argued we have concluded that we have a sub-algebra of $\so^*(8)$, but the only sub-algebra of $\so^*(8)$ with the same dimension as $\spin(2,6)$ is $\so^*(8)$ itself.} 

One can go through the effort of trying to `squash' the spin algebra down into quaternionic matrices, though it is not particularly nice to write down. We know the Lie algebra for $\so^*(8)$ are those matrices satisfying reversion-transpose skew-symmetry:
\begin{equation}
\alg{so}^*(8) = \{a \in \alg{sl}(\mathbb{H}^4) \> \> | \>\> -\widetilde{a}^\intercal = a  \}.
\end{equation}
A generic element of $\so^*(8)$, written as $A = \sum_{i=1}^{28} a_i M_i$, where $a_i$ are the rotation and boost parameters and $M_i$ the quaternionic basis matrices, will look like
\vspace{0.7em}
\begin{equation*}
\hspace{-2.6em}
    \frac{1}{2}\begin{pmatrix}
     a_1 j & a_3 i+a_4 j+a_5 k+a_2 & a_7 i+a_8 j+a_9 k+a_6 & a_{11} i+a_{12} j+a_{13} k+a_{10} \\
     -a_3 i+a_4 j-a_5 k-a_2 & a_{14} j & a_{16} i+a_{17} j+a_{18} k+a_{15} & a_{20} i+a_{21} j+a_{22} k+a_{19} \\
     -a_7 i+a_8 j-a_9 k-a_6 & -a_{16} i+a_{17} j-a_{18} k-a_{15} & a_{23} j & a_{25} i+a_{26} j+a_{27} k+a_{24} \\
     -a_{11} i+a_{12} j-a_{13} k-a_{10} & -a_{20} i+a_{21} j-a_{22} k-a_{19} & -a_{25} i+a_{26} j-a_{27} k-a_{24} & a_{28} j \\
    \end{pmatrix}.
\end{equation*}
\vspace{0.4em}

\noindent
We can then recast this generic basis into precisely the same basis $L_{ij}$ above, acting explicitly upon $\mathbb{H}^4$, via the following assignments of arbitrary parameters $a\rightarrow \theta$, showing the equivalence between the generic Lie algebra elements $A(\vec{a})$ and $s(\vec{\theta})$:
\begin{equation*}
\begin{split}
a_1&\to \theta _{0,1}-\theta _{2,3}+\theta _{4,5}+\theta _{6,7}, \qquad a_{14}\to \theta _{0,1}-\theta _{2,3}-\theta _{4,5}-\theta _{6,7}, \\
a_{23}&\to \theta _{0,1}+\theta _{2,3}+\theta _{4,5}-\theta _{6,7}, \qquad a_{28}\to \theta _{0,1}+\theta _{2,3}-\theta _{4,5}+\theta _{6,7}, 
\end{split}
\end{equation*}
\vspace{-0.7em}
\begin{equation}
\begin{split}
a_2&\to \theta _{4,6}-\theta _{5,7}, \qquad a_3\to \theta _{0,3}-\theta _{1,2}, \qquad \>\>\>\>\>a_4\to \theta _{5,6}-\theta _{4,7},\\
a_5&\to \theta _{1,3}-\theta _{0,2}, \qquad a_6\to \theta _{2,6}-\theta _{3,7}, \qquad \>\>\>\>\>a_7\to \theta _{1,4}-\theta _{0,5},\\
a_8&\to \theta _{3,6}-\theta _{2,7}, \qquad a_9\to \theta _{1,5}-\theta _{0,4}, \qquad \>\>\>\>\>a_{10}\to \theta _{3,5}-\theta _{2,4}, \\
a_{11}&\to \theta _{1,6}-\theta _{0,7}, \qquad a_{12}\to \theta _{2,5}-\theta _{3,4}, \qquad\>\>\> a_{13}\to \theta _{1,7}-\theta _{0,6},\\
a_{15}&\to \theta _{2,4}+\theta _{3,5}, \qquad a_{16}\to -\theta _{0,7}-\theta _{1,6}, \qquad a_{17}\to \theta _{2,5}+\theta _{3,4}, \\
a_{18}&\to \theta _{0,6}+\theta _{1,7}, \qquad a_{19}\to \theta _{2,6}+\theta _{3,7}, \qquad \>\>\>\>a_{20}\to \theta _{0,5}+\theta _{1,4}, \\
a_{21}&\to \theta _{2,7}+\theta _{3,6}, \qquad a_{22}\to -\theta _{0,4}-\theta _{1,5}, \qquad a_{24}\to -\theta _{4,6}-\theta _{5,7}, \\
a_{25}&\to \theta _{0,3}+\theta _{1,2}, \qquad a_{26}\to -\theta _{4,7}-\theta _{5,6}, \qquad a_{27}\to \theta _{0,2}+\theta _{1,3}.
\end{split}
\end{equation}
This gives an explicit construction identifying $\spin(2,6) \cong \so^*(8)$. In this case, we do not have a double covering at the level of the groups, simply isomorphism $\grp{SO}^*(8) \cong \grp{SemiSpin}(2,6)$, because we can see the algebra representations are identified, and so the group representations will be as well via the exponential map. It is worth mentioning again that just as in the Euclidean case, the $L$ and $R$ representations are formally \textit{semi-spin} representations: though they represent connected irreducible double coverings of $\grp{SO}(2,6)$, both only have a center of $\mathbb{Z}_2 = \pm 1$, and so are not faithful representations the group $\grp{Spin}(2,6)$. 

\subsection{Triality of the Quaternionic Orthogonal Group \texorpdfstring{$\grp{SO}^*(8)$}{SO*(8)}}
From the above we have bases for the two spin representations $L$ and $R$ of $\spin(2,6) \cong \so^*(8)$. These quaternionic spin representations exhibit the automorphism of triality with the vector representation of $\alg{so}(2,6)$. To make the triality map explicit between these representations, we will firstly apply a change of basis: this is both to make the triality map simpler, and secondly so that the vector representation (the defining representation of $\grp{SO}(2,6)$) is manifestly real. With
\begin{equation}
    U = \grp{diag}\{i, i, 1,1,1,1,1,1\}, \qquad P = \grp{diag}\{-1, 1, 1,1,1,1,1,1\},
\end{equation}
the changes of basis are
\begin{equation}
    L_{ij} \mapsto U^{-1} L_{ij}\> U, \qquad R_{ij} \mapsto (UP)^{-1} R_{ij} \>U P,
\end{equation}
The change of basis $P$ shows up in the construction of Euclidean and and Lorentzian triality as well. Its necessity can be understood as follows: in the vector case $P$ would be the aforementioned reflection of the $0^{\textrm{th}}$ coordinate relating $L$ and $R$, and as the triality will `hold onto' this change of basis while mapping between representations, without $P$, triality would map $L$ and $R$ to different bases of the vector representation, related by this reflection. As such we can make a choice to preemptively absorb it, only slightly twisting the nice relationship between the $L$ and $R$ bases. Furthermore, the change of basis $U$ appears explicitly to ensure the vector representation is real. It explicitly transforms the bases of $L$ and $R$ into ones which respect the transformed orthogonal structure $U^\intercal g\> U = \mathbb{I}_{2,6}$. 

The triality map in this case is not quite so homogeneous as in the Euclidean or Lorentzian cases. In the bases given, the triality may be explicitly realized in the following way. Collect the $28$ basis elements of these representations into six families `$B$', of four commuting generators; two boosts and two rotations each, and one family of the remaining four compact generators, `$B^\prime$':
\begin{equation}
\label{TrialityQuartets2.0}
B_L=\begin{pmatrix}
    \vec{a} \\
\vec{b} \\
\vec{c} \\
\vec{d}
\end{pmatrix}_{\hspace{-0.5em}L},\quad \mathrm{and} \quad     B_L^\prime = \begin{pmatrix}
    L_{01} \\
    L_{23}\\
    L_{45}\\
    L_{67}\\
\end{pmatrix},\quad \mathrm{with} \qquad
\begin{matrix}
 \vec{a}_L = \{L_{02},L_{03},L_{04},L_{05},L_{06},L_{07}\}, \\
\vec{b}_L = \{L_{13},L_{12},L_{15},L_{14},L_{17},L_{16}\},\\
\vec{c}_L = \{L_{57},L_{47},L_{37},L_{36},L_{24},L_{25}\},\\
\vec{d}_L = \{L_{46},L_{56},L_{26},L_{27},L_{35},L_{34}\},
\end{matrix}
\end{equation}
with $B_R$, $B_V$, $B^\prime_R$, and $B^\prime_V$ defined analogously. Now define matrices which will act on these quartets:
\begin{equation}
H=\frac{1}{2}
    \begin{pmatrix}
    -1 & -1 & 1 & 1 \\
    1 & 1 & 1 & 1 \\
    -1 & 1 & 1 & -1 \\
    -1 & 1 & -1 & 1 \\
    \end{pmatrix},
\qquad G=\frac{1}{2}
\begin{pmatrix}
-1 & -1 & i & -i \\
 1 & 1 & i & -i \\
 i & -i & 1 & 1 \\
 -i & i & 1 & 1 \\
\end{pmatrix}.
\end{equation}
Then the following maps applied simultaneously on this $28$ dimensional space of generators is a Lie algebra automorphism, which cyclically permutes the $V, L, R$ representations of $\spin(2,6)$:
\begin{equation}
B^\prime_{V/L/R} = H\> B^\prime_{L/R/V}, \quad \mathrm{and}\quad B_{V/L/R} = G \>B_{L/R/V}.
\end{equation}
Of note is that $H$ is precisely the triality matrix of the Euclidean case, and appears here again to transform the only quartet containing wholly compact generators. $G$ however is new, as it must transform two rotations and two boosts, into a reordering of the same. In either case the maps cube to the identity, and precisely map our constructed bases of these three representations to one another. For completeness, an arbitrary element of the vector representation which falls out of this triality map is given below.
\begin{equation}
\sum_{0\leq i < j \leq 7} \theta_{i,j}V_{i,j} = \begin{pmatrix}
 0 & \theta _{0,1} & \beta _{0,2} & -\beta _{0,3} & \beta _{0,4} & \beta _{0,5} & \beta _{0,6} & \beta_{0,7} \\
 -\theta _{0,1} & 0 & -\beta _{1,2} & -\beta _{1,3} & \beta _{1,4} & -\beta _{1,5} & \beta _{1,6} & -\beta _{1,7} \\
 \beta _{0,2} & -\beta _{1,2} & 0 & \theta _{2,3} & \theta _{2,4} & -\theta _{2,5} & \theta _{2,6} & -\theta _{2,7} \\
 -\beta _{0,3} & -\beta _{1,3} & -\theta _{2,3} & 0 & \theta _{3,4} & \theta _{3,5} & \theta _{3,6} & \theta _{3,7} \\
 \beta _{0,4} & \beta _{1,4} & -\theta _{2,4} & -\theta _{3,4} & 0 & -\theta _{4,5} & -\theta _{4,6} & \theta _{4,7} \\
 \beta _{0,5} & -\beta _{1,5} & \theta _{2,5} & -\theta _{3,5} & \theta _{4,5} & 0 & -\theta _{5,6} & -\theta _{5,7} \\
 \beta _{0,6} & \beta _{1,6} & -\theta _{2,6} & -\theta _{3,6} & \theta _{4,6} & \theta _{5,6} & 0 & -\theta _{6,7} \\
 \beta_{0,7} & -\beta _{1,7} & \theta _{2,7} & -\theta _{3,7} & -\theta _{4,7} & \theta _{5,7} & \theta _{6,7} & 0 \\
\end{pmatrix}.
\end{equation}
Here non-compact parameters have been superfluously renamed to $\beta$ for clarity, and by inspection one can see this basis obeys the expected orthogonal structure $\mathbb{I}_{2,6} \left(-V^\intercal_{ij}\right) \mathbb{I}_{2,6} = V_{ij}, \>\> \forall i,j$. The three sets of basis elements for the three representations $L_{ij}$, $R_{ij}$, and $V_{ij}$ all obey identical commutation relations. 

Some final comments: just as in \cite{Triality}, if one attempts to find a Lie algebra basis where the triality map acts diagonally, these changes of basis are complex and so will modify the signature of the underlying orthogonal structure, leaving us with a basis naturally understood as spanning $\spin(1,7)$, and the Lorentzian case remains the more `natural' setting for exploring the triality; this is likely related to the relationship between triality and the octonions, which have the same signature. One may also wonder how it is that merely `reordering' the bases, as done via the triality maps, causes the quaternionic structures upon the $L$ and $R$ representations to `become' a real structure upon $V$, and especially what this could mean for the generators when written explicitly in the quaternionic setting. I do not have a complete and compelling answer to this, if there even is one. The main obstruction to an analysis of this kind is how one might attempt to move the various necessary matrices $H, G, U,$ and $P$ into the quaternionic setting. $H$ is no problem, as it is strictly real. $U$ is tantamount to a change of basis making the default reversion structure $\mathbb{I}_{3,1}$, instead of the identity. The other two matrices are more problematic: for $G$, which imaginary unit should $i$ be mapped to? Perhaps if one looks instead at the bi-quaternions $\mathbb{H}\otimes \mathbb{C}$, there is an obvious choice of embedding compatible with everything above (the $i$ in the complex setting mapped to the new commuting imaginary unit), but even then it is not clear how one could map the change of basis $P$ to the quaternionic setting. 

\noindent\rule{16.5cm}{0.4pt}
\appendix
\section{Quaternionic Matrices and Representations}
\label{app:A}
Recall the definition of a quaternion $q \in \mathbb{H}$:
\begin{equation}
    q = t1 + x i+yj +zk, \qquad \textrm{with} \qquad i^2=j^2 = k^2 = ijk = -1 \qquad \textrm{and} \qquad \{t,x,y,z\}\in \mathbb{R}.
\end{equation}
It can be seen from the definition the imaginary units mutually anti-commute: $ij = -ji$ etc. The non-commutativity of quaternions means formally they are not a \textit{field} of numbers, as $\mathbb{R}$ or $\mathbb{C}$ are, but a \textit{normed division ring}. Due to this non-commutativity, many of our favorite facts taken for granted in real and complex representation theory are no longer true, and must be either abandoned or altered to work with quaternions. We will see below that we lose the existence of both the conjugate representation, and the contragredient (or dual) representation, since conjugation and transposition are no longer (anti)homomorphisms on quaternionic matrix algebras.
\subsection{Homomorphisms and Anti-homomorphisms}
A homomorphism from a group $G$ to a group $M$ is a map $T: G\mapsto M$, that commutes with the group multiplications. Specifically if
\begin{equation}
    T(gh) = T(g)T(h) ,\qquad \forall g,h \in G ,
\end{equation}
then $T$ is a group homomorphism. The definition is similar for an algebra. In particular for a Lie algebra $\alg{g}$ and some map $T$ to another algebra $\alg{m}$, if
\begin{equation}
    T([\gamma, \eta]) =[T(\gamma), T(\eta)] , \qquad \forall \gamma, \eta \in \alg{g},
\end{equation}
then $T$ is an algebra homomorphism. An automorphism is a homomorphism with trivial kernel (I.e. $T(g) \in M$ is the identity in $M$ if and only if $g$ is the identity in $G$, and analogously for algebras).  

Homomorphisms are important in representation theory twice over. Firstly because a representation of a group (algebra) upon a vector space, is by definition a homomorphism into the general linear group $\grp{GL}(V,\mathbb{K})$ (endomorphism ring) of that vector space. Secondly they are important as any generic homomorphisms of our matrix rings will in principle induce additional representations. Matrices over the complex numbers have three non-trivial homomorphisms that come `for free': the inverse-transpose map, complex-conjugation, and their composition. This means for any representation we define over a complex vector space, we immediately also have the dual, conjugate, and conjugate dual representations as well (though very often some or all of these coincide with one another). Thus understanding the relevant homomorphisms for quaternionic matrices will help us in understanding quaternionic groups and representations thereof.

Another important concept is the anti-homomorphism, which is like a homomorphism except it reverses the order of multiplication:
\begin{equation}
    T(gh) = T(h)T(g) ,\qquad \forall g,h \in G,
\end{equation}
and analogously for algebras. For matrix groups over $\mathbb{C}$ and $\mathbb{R}$, transposition is an anti-homomorphism, as is inversion. When the image of $T$ is an Abelian group, an anti-homomorphism is of course the same as a homomorphism. Finally the composition of an even number of anti-homomorphisms will be a homomorphism, and so they are as interesting to us as homomorphisms. 

\subsubsection{Conjugation}
Recall quaternionic conjugation of a quaternion $q= t + x i+yj +zk$ is defined as:
\begin{equation}
     q^* = t - xi -yj -zk.
\end{equation}
Conjugation over $\mathbb{H}$ is an anti-homomorphism, as opposed to the case of $\mathbb{C}$ where it is a homomorphism. Were we to assume it is merely a homomorphism we can see immediately a contradiction:
\begin{equation}
    -i = i^* = (jk)^* \stackrel{!}{=} (j^*)(k^*) = (-j)(-k) = jk = i.
\end{equation}
The problematic step where we assume it is a homomorphism is marked with the exclamation, and the contradiction is resolved by instead having conjugation act as an anti-homomorphism 
\begin{equation}
    (jk)^* = (k^*)(j^*).
\end{equation}
Conjugation of a quaternion, when mapped via the standard identification in the complex setting Eq.~(\ref{eq:quat2Pauli}), corresponds to the conjugate transpose of the corresponding matrix.

\subsubsection{Reversion and the Self-Conjugacy of Quaternion Representations}
\label{sec:Reversion}
Should one wish to conjugate only a single imaginary unit, such as the $j$ component, this can be accomplished as $-j q^* j$, and analogously for $i, k$. It may seem overkill to give this operation a name, since it is linearly related to standard quaternion conjugation, but it will prove useful to have language for it. The \textbf{reversion-conjugate} of a quaternion $q=t + x i+yj +zk$ is defined as:\footnote{One reference for the introduction of reversion is \href{https://arxiv.org/pdf/math/0106063}{Introduction to Arithmetic Groups}, by Dave Witte Morris, in particular page 430, example A2.4. \cite{ArithGroups}}
\begin{equation}
\widetilde{q} := t + xi -yj +zk.
\end{equation}
As a composition of maps we have that: $\widetilde{\>\>\>} = j \>\circ\> * \>\circ\> j $. When one considers quaternions as a Cayley-Dickson construction, i.e. as pairs of two complex numbers: $(a+ib) + (c+id)j$, reversion may be thought of as the `conjugation' of the second complex number. It is clear that if quaternion conjugation is an anti-homomorphism, then so too is reversion. 

As noted above we could have chosen any unit for this `additional' (but redundant) conjugation; the choice of $j$ is related to the fact that when one makes the correspondence between $\grp{M}_n (\mathbb{H}) \cong M_{2n}(\mathbb{C})$ via the standard identification in Eq.~(\ref{eq:quat2Pauli}), the imaginary unit $j$ is mapped to the anti-symmetric component of the matrix, while the other components lay in the symmetric subspace. As such reversion of a quaternion corresponds with transposition of its corresponding complex matrix. With standard quaternion conjugation corresponding to the conjugate-transpose of the corresponding complex matrix, the composition of the two: $\widetilde{q}^* = -j q j$ corresponds to \textit{complex} conjugation of the corresponding matrix. For any matrix $A$ defined over the quaternions, and its corresponding matrix mapped into the complex setting $M$, it must be true by definition that
\begin{equation}
    \widetilde{A}^* = -j A j, \qquad \Rightarrow \qquad M^* = JMJ^{-1}.
\end{equation}
Thus any group or algebra defined over quaternionic matrices, when embedded as complex matrices, must be linearly related to its complex conjugate: all representations of quaternionic matrix groups will be self-conjugate. Said another way the heart of why quaternion representations are self-conjugate is that complex conjugation in the complex setting, is a linear map in the quaternionic setting, so any group defined over the quaternions will find its representations are self-conjugate. The relationship to the existence of a quaternionic anti-linear structure is simply a multiplication by $j$ away, since here the operation $j\>\widetilde{()}^*$ is trivially an intertwiner of any group or algebra, since by definition $j \>\widetilde{A}^* =  A \>j$, and the square of this operator is $-1$.

The fact that the composite operation of reversion-conjugation in the quaternionic setting performs complex conjugation in the complex setting, holds even when considering the standard identification of quaternionic and complex vector spaces. Upon a quaternion valued vector $\psi$, one of course has the same fact $ \>\widetilde{\psi}^* =  -j\psi j$. However if we consider the vector space strictly as a left quaternion module, then reversion conjugation is a non-trivial homomorphism we have access to, since we are not able to have scalars multiplying on the right of vectors. Inspecting this action on $\mathbb{H}$ and sending the result to $\mathbb{C}^2$ :
\begin{equation}
    \begin{split}
        \psi &= t+ix+jy+kz \quad \mapsto \quad \begin{pmatrix}t+i \\
    y+iz\end{pmatrix},\\
    \widetilde{\psi}^* &= \widetilde{(t+ix+jy+kz)}^* = t-ix+jy-kz \quad \mapsto\quad  \begin{pmatrix}t-ix \\
    y-iz\end{pmatrix}.\\
    \end{split}
\end{equation}
Not only does conjugation-reversion of quaternion matrices correspond to complex conjugation of the corresponding complex matrices, conjugation-reversion of a quaternion vector can also be chosen to correspond to complex conjugation of the corresponding complex vector. Again this map is no different than the application of $-j()j$, but if one has reason to recast all expressions into a left quaternion module (which is often useful when mapping to the real and complex settings), then reversion an extremely helpful concept for recasting anti-linear operations into the quaternionic setting. 

\subsubsection{Conjugation and Transposition Do Not Play Nice with Quaternion Matrices}
Consider two generic square quaternionic matrices $A$ and $W$. It can be seen that conjugation are transposition separately are neither homomorphisms or anti-homomorphisms: $(A W)^* \neq A^* W^* \neq W^* A^*$ and $(A W)^\intercal \neq A^\intercal W^\intercal \neq W^\intercal A^\intercal$. For the sake of demonstration take a simple example for $M_2(\mathbb{H})$ (recall each entry is a quaternion):
\begin{equation}
\label{eq:2by2}
\begin{pmatrix}
a & b \\
c & d 
\end{pmatrix} 
\begin{pmatrix}
w & x \\
y & z 
\end{pmatrix} = \begin{pmatrix}
aw + by & ax + bz \\
cw + dy & cx + dz 
\end{pmatrix}.
\end{equation}
\subsubsection*{Conjugation}
Quaternionic conjugation of say, the first element on the right hand side would give us $w^* a^* + y^* b^*$. As this flips the order of the component multiplication, for the left hand side to have a chance of equaling the conjugated right side, we would need to extend conjugation on quaternion matrices as an anti-homomorphism. Applying this assumption to each side of the equation gives:
\begin{equation}
\begin{split}
\begin{pmatrix}
w^* & x^* \\
y^* & z^* 
\end{pmatrix}
\begin{pmatrix}
a^* & b^* \\
c^* & d^* 
\end{pmatrix} &\stackrel{?}{=} \begin{pmatrix}
w^*a^* + y^*b^* & x^*a^* + z^*b^* \\
w^*c^* + y^*d^* & x^*c^* + z^*d^* 
\end{pmatrix},\\ 
\begin{pmatrix}
w^*a^*+ x^*c^* & w^*b^* + x^*d^* \\
y^*a^* + z^*c^* & y^*b^* + z^*d^* 
\end{pmatrix} &\stackrel{?}{=} \begin{pmatrix}
w^*a^* + y^*b^* & x^*a^* + z^*b^* \\
w^*c^* + y^*d^* & x^*c^* + z^*d^* 
\end{pmatrix}.\\ 
\end{split}
\end{equation}
This is not true in general, even if it looks like the correct answer merely shuffled around. Thus although quaternion conjugation is an anti-homomorphism of $\mathbb{H}$, it cannot be made to be a homomorphism nor an anti-homomorphism on $\grp{M}_n(\mathbb{H})$. The same results follow for reversion, since it is linearly related to quaternion conjugation. This means for a generic representation of a quaternion valued matrix group or algebra, conjugation is not a homomorphism and so \textbf{there are no quaternion-conjugate representations}.

\subsubsection*{Transposition}
Looking again at the simple example Eq.~(\ref{eq:2by2}), if we assume transposition is an anti-homomorphism as usual, on the left we would necessarily reverse the order in which the matrix elements are being multiplied, while on the right we would not. So if transposition is going to stand a chance of having any nice properties over quaternionic matrices, we should assume it acts as a homomorphism. Applying it to each side of the equation we have
\begin{equation}
\begin{split}
\begin{pmatrix}
a & c \\
b & d 
\end{pmatrix}
\begin{pmatrix}
w & y \\
x & z 
\end{pmatrix} &\stackrel{?}{=} \begin{pmatrix}
aw + by & cw + dy \\
ax + bz & cx + dz 
\end{pmatrix}, \\ 
\begin{pmatrix}
aw+cx & ay + cz \\
bw + dx & by + dz 
\end{pmatrix} &\stackrel{?}{=} \begin{pmatrix}
aw + by & cw + dy \\
ax + bz & cx + dz 
\end{pmatrix}.
\end{split}
\end{equation}
This is not true in general, and so transposition cannot be an (anti)homomorphism over quaternion matrices. This means for a generic representation of a quaternion valued matrix group or algebra, we cannot invoke the standard inverse-transpose homomorphism and so \textbf{there is no contragredient (dual) representation}.

\subsubsection{The Conjugate Transpose is an Anti-Homomorphism}
The above failures to construct standard matrix (anti)automorphisms come from the non-commuting nature of the quaternions. Inspection of just how close the above actions are to being satisfactory begs one to compose them --- and indeed the composition of the above failures becomes a success!\footnote{Never let anyone tell you two wrongs cannot make a right.} Inspecting the example above once again, taking the quaternion Hermitian conjugate of both sides:
\begin{equation}
\begin{split}
\begin{pmatrix}
w^* & y^* \\
x^* & z^* 
\end{pmatrix}
\begin{pmatrix}
a^* & c^* \\
b^* & d^* 
\end{pmatrix} &\stackrel{?}{=} \begin{pmatrix}
w^*a^* + y^*b^* & w^*c^* + y^*d^* \\
x^*a^* + z^*b^* & x^*c^* + z^*d^* 
\end{pmatrix},\\ 
\begin{pmatrix}
w^*a^* + y^*b^* & w^*c^* + y^*d^* \\
x^*a^* + z^*b^* & x^*c^* + z^*d^* 
\end{pmatrix} &= \begin{pmatrix}
w^*a^* + y^*b^* & x^*a^* + z^*b^* \\
w^*c^* + y^*d^* & x^*c^* + z^*d^* 
\end{pmatrix}.
\end{split}
\end{equation}
Thus even without dual or conjugate representations, we still have the conjugate-dual representation. For quaternionic matrices we have just shown: 
\begin{equation}
    (AW)^\dagger = W^\dagger A^\dagger,
\end{equation}
which of course implies the reversion-transpose is an anti-homomorphism as well
\begin{equation}
     \widetilde{(AW)}^\intercal = -j(AW)^\dagger j = -j W^\dagger A^\dagger j = (-j W^\dagger j) (-jA^\dagger j) = \widetilde{W}^\intercal \widetilde{A}^\intercal.
\end{equation}
The conjugate dual representation will exist for algebras too. If we have some real Lie algebra with a representation over quaternionic matrices, with a bracket in some basis:
\begin{equation}
    \left[\lambda_i, \lambda_j\right] = f^{ijk} \lambda_k.
\end{equation}
Then we can see:
\begin{equation}
    \begin{split}
        \left[-\lambda_i^\dagger, -\lambda_j^\dagger\right] &= (\lambda_i^\dagger \lambda_j^\dagger - \lambda_j^\dagger \lambda_i^\dagger) = -(\lambda_i \lambda_j - \lambda_j\lambda_i)^\dagger \\
        & = -\left[\lambda_i, \lambda_j\right]^\dagger =f^{ijk} (-\lambda_k^\dagger),
    \end{split}
\end{equation}
the action $-\lambda^\dagger$ is an algebra homomorphism, giving a basis for the conjugate-dual representation.

\noindent\rule{16.5cm}{0.4pt}
\section{Groups over the Quaternions}
\label{app:B}
There are four types of quaternionic matrix groups. They are named:
\begin{equation}
    \grp{GL}(\mathbb{H}^n), \qquad \grp{SL}(\mathbb{H}^n), \qquad \grp{Sp}^*(p,q, \mathbb{H}^{p+q}), \qquad \grp{SO}^*(2n,\mathbb{H}^n).
\end{equation} 
The naming conventions for the latter two are rather confusing when in the arena of quaternionic matrices; the names refer to the behaviour of the groups \textit{after} complexification.\footnote{Complexification of a quaternionic matrix algebra is done first by mapping to $\grp{M}_{2n}(\mathbb{C})$ via Eq.~(\ref{eq:quat2Pauli}), and then allowing complex combinations of the Lie algebra generators.} Synonyms for the groups might be `Unitary-quaternion group' for $\grp{Sp}^*$, and `Skew-unitary' or `reversion unitary' for $\grp{SO}^*$. We will define these groups below.

\subsection{\texorpdfstring{$\grp{GL}(\mathbb{H}^n)$}{GL(n, H)}}
The General Linear Group on quaternions is no problem to define, it is simply every linear invertible quaternion valued matrix acting on some $\mathbb{H}^n$, acting from the left, in terms of both the matrix and quaternion multiplication. Simple counting gives a real dimension of $4n^2$.

\subsection{\texorpdfstring{$\grp{SL}(\mathbb{H}^n)$}{SL(n, H)}}
The Special Linear Group on quaternions is not so easy to define. Because the elements generally do not commute there are as many ways to define a determinant as elements multiplied in said determinant, which grows with the factorial of the dimension. The standard approach to deal with this is known as the `Study Determinant': simply send the quaternionic matrix to its complex counter part in twice the dimension, and take the determinant there.\footnote{For more on quaternion matrix determinants, see the \href{https://ncatlab.org/nlab/show/Dieudonn\%C3\%A9+determinant}{Dieudonné determinant} \cite{det} on nlab, and citations therein.} The corresponding Lie algebra $\alg{sl}(\mathbb{H}^n)$ is given by those quaternionic matrices which have a purely imaginary trace --- this corresponds to traceless-ness in the complex case (recall the embedding Eq.~(\ref{eq:quat2Pauli})). This gives a real dimension of $4n^2 -1$. 

In the interest of keeping track of real forms we can gather more information. Inspection of an arbitrary traceless quaternion matrix reveals its largest compact sub-algebra will be of real dimension $n(2n+1)$, which we will see later is the group $\grp{Sp}^*(n, \mathbb{H})$. This implies the algebra $\alg{sl}_n(\mathbb{H})$ has a Killing form $K$ with signature:
\begin{equation}
    \left[n_-, n_+\right] = \left[n(2n+1), \>(n-1)(2n+1)\right], \qquad \operatorname{Tr}(K) = -2n-1.
\end{equation}

\subsection{Every Bi-linear Map on Quaternions is 0}
The remaining families of quaternionic matrix groups are defined below, but first it is fun to consider why there are not the standard generalizations of the \textit{bi-linear} orthogonal and symplectic groups, which exist in both the real and complex cases. If the non-existence of a dual representation is not enough to convince one that quaternionic groups preserving a bi-linear form cannot exist, then the following argument will. Suppose for some $n$ dimensional space of quaternions $\mathbb{H}^n$,  we are looking to preserve a symmetric (orthogonal) bi-linear form $G: \mathbb{H}^n\times\mathbb{H}^n \rightarrow \mathbb{H}$. Starting with vectors $a,b$ (column vectors with quaternion entries), by bi-linear we mean:
\begin{equation}
G\left(i a,j b\right) = ijG\left(a,b\right)
\end{equation}
since $i,j$ are scalars. There are other possible definitions of bi-linearity, but no matter how you try to demand it you will find the same result. Utilizing the symmetry of the form, we have:
\begin{equation}
G\left(i a,j b\right) = G\left(jb,ia\right)= jiG\left(b,a\right)=-ij G\left(a,b\right)
\end{equation}
in the last equality we anti-commute $i$ and $j$, and swap the arguments. Thus we have shown
\begin{equation}
ij G\left(a,b\right)=-ijG\left(a,b\right), \qquad \Rightarrow \qquad 
G\left(a,b\right)=-G\left(a,b\right).
\end{equation}
This tells us that the only symmetric bi-linear form over the quaternions is $0$. The argument is identical for symplectic (anti-symmetric) bi-linear forms since we swap the arguments twice. As such any bi-linear form over quaternions, in any dimension, is identically zero. One cannot even begin to build something like a classical orthogonal or symplectic group over the quaternions. 

\subsection{\texorpdfstring{$\grp{Sp}^*(p,q, \mathbb{H})$}{Sp*(p,q, H)}: Preserving Quaternionic Hermitian Forms}
The quaternionic unitary groups $\grp{Sp}^*(p,q)$ in their defining representation are matrices which preserve a Hermitian sesqui-linear inner product $H: \mathbb{H}^n\times\mathbb{H}^n \rightarrow \mathbb{H}$, where the first argument of $H(a,b)$ is quaternion conjugated, i.e. $H(a,b) = a^{\dagger}hb$ where $h$ is the matrix of $H$ and $q^\dagger = (q^\intercal)^*$. Let's assume we are in a standard basis where $h = \mathbb{I}_{p,q}$ with $p+q=n$. Then we may define the group of matrices which are `quaternion unitary' with respect to $h$ via:
\begin{equation}
\grp{Sp}^*(p,q) = \grp{U}(p,q,\mathbb{H}) = \{A \in \grp{SL}(n, \mathbb{H}) \> \>| \>\> A^\dagger \mathbb{I}_{p,q} A = \mathbb{I}_{p,q}\}
\end{equation}
The star in the signifies that this group is defined over the quaternions, and the name comes from that fact that when complexified, the group becomes the symplectic group $\grp{Sp}(2n, \mathbb{C})$. When the signature is definite (one of $p$ or $q$ is $0$) the group is compact. 

It is well known that the compact symplectic group can be understood to be the intersection of the unitary and symplectic groups: $\mathrm{U}(2n, \mathbb{C}) \>\cap \> \mathrm{Sp(2n, \mathbb{C})}$. This fact can be seen generally for any signature as follows: by the definition of reversion, if the group elements obey $A^\dagger \mathbb{I}_{p,q} A = \mathbb{I}_{p,q}$, they must also obey $\widetilde{A}^\intercal \left(j \>\mathbb{I}_{p,q}\right) A = j\>\mathbb{I}_{p,q}$. When these conditions are translated to the complex setting of twice the dimension, the first condition corresponds to being a complex (pseudo)unitary matrix in $\grp{U}(p,q)$; the second condition corresponds to being in the complex symplectic group $\grp{SP}(2n, \mathbb{C})$. 

The dimension will be the same for the various signatures, so it is easier to consider the compact case first. The algebra $\alg{sp}^*_n(\mathbb{H})$ will be all anti-Hermitian (Hermitian in the quaternionic sense) quaternionic matrices, which gives $n(n-1)/2$ real anti-symmetric parts, $3 n(n+1)/2$ imaginary symmetric parts, giving a total dimension of $d =n(2n+1)$. The largest compact subgroup for $\grp{Sp}^*(p,q)$ will be given by $\grp{Sp}^*(p) \times \grp{Sp}^*(q)$, which has dimension:
\begin{equation}
    p(2p+1) + q(2q+1) = 2(p^2+q^2)+(p+q) = n(2n+1) - 4pq.
\end{equation}
This means the Killing form $K$ of the corresponding algebra will have signature:
\begin{equation}
    \left[n_-, n_+\right] = \left[n(2n+1) - 4pq, \> 4pq\right], \qquad \operatorname{Tr}(K) = 8pq - n(2n+1).
\end{equation}

\subsection{\texorpdfstring{$\grp{SO}^*(2n, \mathbb{H}^n)$}{SO*(2n, H)}: Quaternionic Orthogonal Groups}
\label{sec:quatorthgrps}
Finally we reach the least well known set of real forms of the classical groups: the so-called `quaternionic orthogonal' groups \cite{book:symrepinv}. Seeing as we recently demonstrated the impossibility of such groups, it seems like a terrible name. Regardless, they are so named because after complexification they become a complex orthogonal group $\grp{SO}(2n, \mathbb{C})$. This group is defined as the preservation of a \textit{reversion} symmetric sesqui-linear inner product $G: \mathbb{H}^n\times\mathbb{H}^n \rightarrow \mathbb{H}$: where the first argument of $G(a,b)$ undergoes a reversion i.e. $G(a,b) = \widetilde{a}^{\intercal}gb$ where $g$ is the matrix of $G$. Let's assume we are in a standard basis where $g = \mathbb{I}_{n}$. Then we may define the group of matrices which are `reversion unitary' with respect to $g$ via:
\begin{equation}
\grp{SO}^*(2n, \mathbb{H}^n) = \{O \in \grp{SL}(n, \mathbb{H}) \>\> | \>\> \widetilde{O}^\intercal \mathbb{I}_n O = \mathbb{I}_{n}\}.
\end{equation}
The dimension of $\grp{SO}^*(2n)$ is the same as an even orthogonal group: canonically the Lie algebra would have $n(n-1)/2$ independent quaternions on the off diagonal, and $n$ quaternions on the diagonal which only have $j$ components, leading to a real dimension of $n(2n-1)$. As for the largest compact sub-algebra, it can be seen that the off diagonal entries attached to $i$ or $k$ components give rise to imaginary-anti-symmetric (Hermitian) generators, and thus a Killing form $K$ of signature: 
\begin{equation}
    \left[n_-, n_+\right] = \left[n^2, n(n-1)\right], \qquad \operatorname{Tr}(K) = -n.
\end{equation}
With all this stated, I imagine a few questions may pop into the readers head. Firstly, how can this be distinct from the quaternion-conjugate groups if reversion is linearly related to standard conjugation? Secondly, are there not forms of this group with signature $p,q$? 

The answer to the first question is given in the next section. The answer to the second question is that Sylvester's law of inertia\footnote{The question of which transformations leave the number of positive and negative eigenvalues (including multiplicity) in tact, for a particular structure. See \cite{book:symrepinv} and \cite{ArithGroups} for review.} no longer holds for the reversion inner product. Since $j$ is the only imaginary unit being conjugated, one can perform a change of basis involving $i$ or $k$ which leaves $g$ with any signature one pleases. This is exactly analogous to how when we bring the real (bi-linear form preserving) orthogonal groups to the complex realm, all choices of signature are a change of basis away from one another. 

\subsection{Skew-Hermitian Forms: The Relationship Between \texorpdfstring{$\grp{SO^*}$}{SO*} and \texorpdfstring{$\grp{Sp}^*$}{Sp*}}
\label{sec:waxdiff}
The key dis-analogy between the complex case and the quaternionic case, is that in the complex case, preserving anti-Hermitian or Hermitian forms lead to the same group, since every anti-Hermitian form is merely a multiplication by $i$ away from being a Hermitian form, and scalar multiplication commutes. Given the defining group relation for some unitary group over the complex numbers:
\begin{equation}
    U^\dagger h U = h, \quad \Rightarrow \quad U^\dagger (ih )U = (ih),
\end{equation}
it is clear the notion of a skew-unitary group is redundant: the definitions define the same groups. However due to the non-commutativity of quaternions, this trick does not work for quaternion-Hermitian forms:
\begin{equation}
    A^\dagger h A = h, \quad \Rightarrow \quad (A^\dagger h A)i = hi,
\end{equation}
the $i$ in general cannot be commuted inside the product. We anticipate then the transformations which preserve quaternion-Hermitian inner products, and transformations which preserve quaternion-skew-Hermitian inner products, in general should define different groups, analogous to how preservation of symmetric and anti-symmetric bi-linear forms lead to two completely distinct groups over the reals. 

To understand the complimentary nature of these two groups, let us inspect how the defining group relations are affected by recasting quaternion conjugation as reversion and vice versa. Let's assume we have a group preserving a quaternion-Hermitian form with a well defined parity: 
\begin{equation}
    h^\dagger = \pm h.
\end{equation} 
Then staring from
\begin{equation}
A^\dagger h A = h, \qquad A\in \grp{SL}_n(\mathbb{H}),
\end{equation} 
we can write the quaternion conjugation as a reversion and multiplication by $j$'s, giving
\begin{equation}
(-j \widetilde{A}^\intercal j) h A = h \qquad \Rightarrow \qquad \widetilde{A}^\intercal (j h) A = jh.
\end{equation}
The group of matrices $A$ which preserved the (skew)Hermitian form $h$ are the same matrices which preserve a `reversion-form' $jh$. It would be particularly nice if the quaternion valued matrix $jh$ had well defined behaviour under reversion-transpose, let's inspect:
\begin{equation}
\left(\widetilde{jh}\right)^\intercal = \widetilde{h}^\intercal\widetilde{j} = - \widetilde{h}^\intercal j = -(-jh^\dagger j)j = -jh^\dagger = \mp jh
\end{equation}
Where in the fourth expression we relate reversion to quaternion conjugation again. If $h$ has well defined symmetry under (quaternion)Hermitian-conjugation, then $jh$ has a well defined and opposite symmetry under reversion-transpose, and the same matrices preserve these two forms simultaneously. 

The takeaway is that each of our two groups simultaneously preserve two structures: $\grp{Sp}^*(p,q,\mathbb{H})$ may equally well be regarded as either the group preserving a standard quaternion-Hermitian form, or equivalently as the group preserving the standard skew-reversion form. These structures respect Sylvester's law of inertia and so each signature is meaningful. Conversely $\grp{SO}^*(2n, \mathbb{H})$ may be interpreted as the group preserving reversion-symmetric quaternion matrices, or equivalently as the group preserving skew-Hermitian quaternion inner products. These structures will not respect Sylvester's law of inertia, and as such every signature is equivalent to the identity. As previously discussed this makes clear the rise of quaternionic structures when we send ourselves from $\grp{M}_n(\mathbb{H})$ to $\grp{M}_{2n}(\mathbb{C})$. In the complex setting, for $\grp{Sp}^*$ one preserves simultaneously a Hermitian form and a symplectic form, while for $\grp{SO}^*$ one preserves simultaneously an anti-Hermitian form and a symmetric bi-linear form. In both cases by compatible triples, this leads to the existence of a quaternionic structure. 

This fact gives another characterization of the quaternion orthogonal groups: they are the unique subgroup (up to change of basis) of the complex orthogonal group $\grp{SO}(2n,\mathbb{C})$ which commute with a chosen quaternionic structure:
\begin{equation}
    \grp{SO}^*(2n) = \{A \in \grp{SO}(2n,\mathbb{C}) \> \>| \>\> J A^* = A J, \quad J J^* = -\mathbb{I}_{2n}\},
\end{equation}
where $^*$ here is entry-wise complex conjugation.


\newpage
\subsection{Table of Quaternion Groups}
\vspace{1ex}
\begin{table}[h]
    \hspace{-4.5em}
    \Centering
    \large
    \begin{tabular}{|c|cccc|}
        \hline Name & dim $d$ & rank & $(n_-, n_+)$ & Complexification\\\hline
        $\alg{sl}(\mathbb{H}^n)$ & $4n^2-1$ & $2n-1$ & $\left[n(2n+1),\>(n-1)(2n+1)\right]$& $\alg{sl}(2n, \mathbb{C})$\\
        $\alg{sp}^*(p,q,\mathbb{H}^{n = p+q})$ & $n(2n+1)$ & $n$ & $\left[d - 4 pq,\>4 pq\right]$& $\alg{sp}(2n, \mathbb{C})$\\
        $\alg{so}^*(\mathbb{H}^n)$ & $n(2n-1)$ & $n$ & $\left[n^2,\>n(n-1)\right]$& $\alg{so}(2n, \mathbb{C})$\\\hline
    \end{tabular}
\begin{changemargin}{-1cm}{-1cm}
    \caption{Compiled facts about real Lie algebras over $\mathbb{H}$: the real dimension and rank of the group, as well as the signature $(n_-, n_+)$ of the Killing form.}
    \label{tab:quatgrps}
\end{changemargin}
\end{table}

\noindent\rule{16.5cm}{0.4pt}
\section{Minutiae}
\label{app:Minutiae}
\subsection*{Pauli Matrices}
The standard Pauli matrices are given below:
\begin{equation}
\label{pauli}
\begin{split}
        \sigma_x = \begin{pmatrix}
            0 & 1 \\
            1 & 0 \\
        \end{pmatrix}, \qquad \sigma_y &= \begin{pmatrix}
            0 & -i \\
            i & 0 \\
        \end{pmatrix}, \qquad \sigma_z = \begin{pmatrix}
            1 & 0 \\
            0 & -1 
        \end{pmatrix}.
\end{split}
\end{equation}

\subsection*{Standard Basis for \texorpdfstring{$\alg{su}(3)$}{su(3)}}
\label{sec:su3}
The standard three dimensional representation of Lie algebra $\su(3)$, which obeys $\lambda = -\lambda^\dagger$, is given below:
\begin{equation}
    \begin{split}
        \lambda_1 = \frac{1}{2}\begin{pmatrix}
            0 & i & 0 \\
            i & 0 & 0 \\
            0 & 0 & 0 
        \end{pmatrix}, \qquad \lambda_2 &= \frac{1}{2}\begin{pmatrix}
            0 & -1 & 0 \\
            1 & 0 & 0 \\
            0 & 0 & 0 
        \end{pmatrix}, \qquad \lambda_3 = \frac{1}{2}\begin{pmatrix}
            i & 0 & 0 \\
            0 & -i & 0 \\
            0 & 0 & 0 
        \end{pmatrix},\\
        \lambda_4 = \frac{1}{2}\begin{pmatrix}
            0 & 0 & i \\
            0 & 0 & 0 \\
            i & 0 & 0 
        \end{pmatrix}, \qquad \lambda_5 &= \frac{1}{2}\begin{pmatrix}
            0 & 0 & -1 \\
            0 & 0 & 0 \\
            1 & 0 & 0 
        \end{pmatrix}, \qquad \lambda_8 = \frac{1}{2\sqrt{3}}\begin{pmatrix}
            i & 0 & 0 \\
            0 & i & 0 \\
            0 & 0 & -2i 
        \end{pmatrix},\\
        \lambda_6 = \frac{1}{2}\begin{pmatrix}
            0 & 0 & 0 \\
            0 & 0 & i \\
            0 & i & 0 
        \end{pmatrix}, \qquad \lambda_7 &= \frac{1}{2}\begin{pmatrix}
            0 & 0 & 0 \\
            0 & 0 & -1 \\
            0 & 1 & 0 
        \end{pmatrix}.
    \end{split}
\end{equation}
\subsection*{\texorpdfstring{$\grp{SO}^*(6)$ Change of Basis Matrix $U$}{SO*(6) Change of Basis Matrix U}}
The change of basis matrix used to diagonalize the Hermitian form $h$ which $\grp{SO}^*(6)$ respects is 
\begin{equation}
\label{changeofbasis}
    U = \frac{1}{\sqrt{2}}
\begin{pmatrix}
 0 & 0 & i & 0 & 0 & -i \\
 0 & 0 & 1 & 0 & 0 & 1 \\
 0 & i & 0 & 0 & -i & 0 \\
 0 & 1 & 0 & 0 & 1 & 0 \\
 i & 0 & 0 & -i & 0 & 0 \\
 1 & 0 & 0 & 1 & 0 & 0 \\
\end{pmatrix}.
\end{equation}

\subsection*{Dirac Matrices \texorpdfstring{$\gamma^\mu$}{ }}
\label{eq:dirac}
\begin{equation}
\gamma^0 = \begin{pmatrix}
 \mathbb{I}_2 & 0 \\
 0 & -\mathbb{I}_2 \\
\end{pmatrix}, \quad \gamma^i = \begin{pmatrix}
    0 & \sigma_i\\
    -\sigma_i & 0
\end{pmatrix}, \quad \gamma^5= \begin{pmatrix}
    0 & \mathbb{I}_2 \\
    \mathbb{I}_2 & 0
\end{pmatrix}.
\end{equation}

\subsection*{Clifford Algebra Generators  \texorpdfstring{$g_i$ for $\grp{C}\ell(7,0)$}{g_i for Cl(7,0)}}
\label{app:cl7}
\begin{equation}
    \begin{split}
        g_1 &= -i \>\mathbb{I}_2 \otimes \gamma^1 \gamma^3 = -i\begin{pmatrix}
            \gamma^1 \gamma^3 & 0 \\
            0 & \gamma^1 \gamma^3,
        \end{pmatrix}\\
        g_2 &=\>\>\> i\sigma_z \otimes \gamma^3 , \qquad g_3 = -i \>\mathbb{I}_2 \otimes \gamma^1,\\
        g_4 &= -i \sigma_y \otimes \gamma^1 \gamma^2, \quad g_5 = \sigma_y \otimes \gamma^1 \gamma^5,\\
        g_6 &= \>\>\>i \sigma_x \otimes \gamma^3, \qquad g_7 = -\sigma_y \otimes \gamma^0 \gamma^1.\\
    \end{split}
\end{equation}
\label{}

\newpage
\bibliographystyle{plain} 
\bibliography{Entire} 

\end{changemargin}
\end{document}